\documentclass[12pt]{amsart}
\usepackage{amsmath, amssymb, amscd}
\usepackage{amssymb}\usepackage[utf8]{inputenc}
\usepackage{mathtools}
\usepackage[all]{xy}

\newtheorem{theorem}{Theorem}[section]
\newtheorem{lemma}[theorem]{Lemma}
\newtheorem{corollary}[theorem]{Corollary}
\newtheorem{proposition}[theorem]{Proposition}
\theoremstyle{definition}

\newtheorem{remark}[theorem]{Remark}

\long\def\comment#1{\relax}
\let\congr=\equiv
\DeclareMathOperator{\Int}{Int}
\newcommand{\vP}{v_{\scriptscriptstyle P}}
\newcommand{\vp}{v_{p}}
\newcommand{\vQ}{v_{\scriptscriptstyle Q}}

\newcommand{\rest}{\negthinspace \restriction\negthinspace}
\newcommand{\gal}{\text{\rm Gal}}

\let\natn=\N

\let\intz=\Z
\newcommand{\OK}{{\mathcal{O}}_K}
\newcommand{\OF}{{\mathcal{O}}_F}
\newcommand{\OL}{{\mathcal{O}}_L}
\newcommand{\OZ}{{\mathcal{O}}_Z}
\newcommand{\OKhash}{\OK^\#}
\newcommand{\OKP}{{\mathcal{O}}_{K,P}}
\newcommand{\OZQZ}{{\mathcal{O}}_{Z, Q_Z}}
\newcommand{\OLQ}{{\mathcal{O}}_{L, Q}}
\newcommand{\IntOK}{\Int(\OK)}

\newcommand{\Spec}{{\mathrm{Spec}}}

\newcommand{\GEn}{{\mathrm{GE}}_n}
\newcommand{\GEt}{{\mathrm{GE}}_2}

\newcommand{\Pc}{{\mathcal{P}}(c)}

\newcommand{\DP}{\negmedspace : \negmedspace}

\renewcommand{\P}{\mathbb P}
\newcommand{\PK}{\P_K}
\newcommand{\PL}{\P_L}

\newcommand{\ito}{\overset\sim\to}

\author{Sophie Frisch}
\address{Institut f\"ur Analysis und Zahlentheorie\\
Graz University of Technology\\ Kopernikusgasse 24\\8010 Graz, Austria}
\email{frisch@math.tugraz.at}
\thanks{S.~Frisch is supported by the Austrian Science Fund (FWF):
P35788}

\author{Franz Halter-Koch}
\address{Institut f\"ur Mathematik und wissenschaftliches Rechnen\\
University of Graz \\ Heinrichstra\ss e 36 \\8010 Graz, Austria}

\title[P-adic approximation of algebraic integers]{P-adic 
approximation of algebraic integers and residue
class rings of rings of integer-valued polynomials}

\begin{document}

\begin{abstract}
Let $F\colon K$ be a Galois extension of number fields and $Q$
a prime ideal of $\OF$ lying over the prime $P$ of $\OK$. By
analyzing the $Q$-adic closure of $\OK$ in $\OF$ we characterize 
those rings of integers $\OK$ for which every residue class ring 
of $\IntOK$ modulo a non-zero prime ideal is $\GEt$ (meaning that
every unimodular pair can be transformed to $(1, 0)$ by a series
of elementary transformations).
\end{abstract}

\maketitle

\section{P-adic approximation of algebraic integers}
For an algebraic number field $K$, let $\OK$ be its ring of integers 
and $\P_K$ the set of all maximal ideals of $\OK$. 
For $P \in \P_K$ we denote by $K_P$ the $P$-adic
completion of $K$ and by $\widehat O_P$ its valuation domain. 

Let $L/K$ be a finite extension of algebraic number fields, and
$P \in \P_K$ and $Q \in \P_L$ such that $Q \supset P$, and 
suppose that $K_P \subset L_Q$. The greatest intermediate field 
$Z$ of $L/K$ satisfying  $Z_{Q \cap Z} =K_P$ is called the 
\,{\bf decomposition field}\, of $Q$ over $K$. 

If $L/K$ is a Galois extension and $G = \gal(L/K)$ then the 
decomposition field $Z$ of $Q$ over $K$ is the fixed field of the 
decomposition group $G_Q = \{\sigma \in G \mid \sigma Q=Q\}$,
the local extension $L_Q/K_P$ is Galois, and the restriction 
$\tau \mapsto \tau \rest L$ defines an isomorphism 
$\gal(L_Q/K_P) \ito G_Q$ 
(see \cite[\S6.1.{\bf 3}]{N} or \cite[Def.~2.5.3]{HKalgnum}). 
We identify: $G_Q = \gal(L_Q/K_P) \subset G$.

In general, the existence of the decomposition field is provided 
by the following simple lemma.
%\smallskip

\begin{lemma}
Let $L/K$ be a finite extension of algebraic number 
fields, $P \in \P_K$ and $Q \in \P_L$ such that $Q \supset P$, 
and suppose that $K_P \subset L_Q$. Then  $K_P\cap L$ is the 
decomposition field of \,$Q$ over $K$.
\end{lemma}

\begin{proof}
If $M$ is an intermediate field of $L/K$ such that 
$M_{Q \cap M} = K_P$, then clearly $M \subset L \cap K_P $, and it 
suffices to prove that $(L \cap K_P)_{Q \cap (L \cap K_P)} = K_P$.
But $K \subset L \cap K_P \subset K_P$, and if we build the 
topological closure (in the $Q$-adic topology), we obtain that
$K_P \subset (L\cap K_P)_{Q \cap L \cap K_P} \subset K_P$, and 
thus equality holds.
\end{proof}
%\medskip

\begin{theorem}\label{1}
Let $L/K$ be a finite extension of algebraic number fields, 
$P \in \P_K$, \ $Q \in \P_L$ such that $Q \supset P$, and let 
$Z$ be the decomposition field of $Q$ over $K$. 

Let $Q_Z = Q\cap \OZ$. Then $\OZQZ$ 
% the localization of $\OZ$ at $Q_Z$ 
is the $Q$-adic closure of $\OK$ in $L$.
\end{theorem}

\begin{proof} 
Let $\overline{\OK}$ be the $Q$-adic closure of $\OK$ in $\OL$. 
By definition, $\overline{\OK}  = \widehat O_P \cap L$,
and thus we obtain: 
\[
\xymatrix@=14pt{
L\ar@{-}[d] \ar@{-}[r] &
L_Q\ar@{-}[r]\ar@{-}[d] &
\widehat O_Q \ar@{-}[d]\\
Z\ar@{-}[d] \ar@{-}[r] &K_P = Z_{Q \cap Z} \ar@{-}[r] &
\widehat O_P = \widehat O_{Q \cap Z}\\
K \ar@{-}[ur]\ar@{-}[rr] &&\OK \ar@{-}[u]
} 
\]
$\overline{\OK}  = 
\widehat O_P \cap L = (K_P \cap \widehat O_Q) \cap L  =
(K_P \cap L) \cap (\widehat O_Q \cap L) = 
Z \cap \OLQ = \OZQZ$.
\end{proof}
%\medskip

\begin{corollary}
Let $L/K$ be a finite Galois extension of algebraic number fields, 
$P \in \P_K$ and $Q \in \P_L$ such that $Q \supset P$,
%\smallskip

\begin{enumerate}
\item
$\OK$ is dense in $\OL$ in the $Q$-adic topology if and only if 
\,$P$ splits completely in $L$.
%\smallskip

\item
$\OK$ is relatively closed in $\OL$ in the $Q$-adic topology 
if and only if \,$Q$ is the only maximal ideal of \,$\OL$ lying above $P$.
\end{enumerate}
\end{corollary}

\begin{proof} 
Let $Z$ be the decomposition field of $Q$ over $K$.
Then $r = [Z \DP K]$ is the number of prime ideals of $\OL$ lying 
above $P$ and $\OZ$ is the $Q$-adic closure of $\OK$ in $\OL$.
Consequently,
\[
\OK \ \text{ is dense in } 
\ \OL \,\iff\, \OZ=\OL\, \iff \,Z=L \, \iff \, r = [L\DP K].
\]
\[
\OK \ \text{ is closed in }
\ \OL \,\iff\, \OZ=\OK \,\iff\, Z=K\, \iff \,r =1.\qedhere
\]
\end{proof}
%\medskip

\begin{theorem}\label{intersectionOK}
Let $L/K$ be a finite extension of algebraic number fields, and 
let $\OKhash$ be the intersection of all $Q$-adic closures of \,$\OK$ 
in $L$. Then $\OKhash = \OK$.
\end{theorem}

\begin{proof} 
Let $K'$ be the intersection of all decomposition fields of $Q$, where 
$Q$ ranges through $\PL$. By Theorem~\ref{1}, $\OKhash\subseteq K'$.
We show that $K' = K$.

Assume first that $L/K$ is Galois and let $G = \gal(L/K)$. 
If $P \in \P_K$ is unramified in $L$ and $Q \in \P_L$ such that 
$Q \supset P$, then $G_Q = \langle F_Q \rangle$ is a cyclic group 
generated by the Frobenius automorphism $F_Q$ of $Q$ over $K$.
Therefore the fixed field $L^{\langle F_Q\rangle}$ is the 
decomposition field of $Q$ over $K$.

For $\sigma \in G$ we denote by $\mathcal P(L/K,\sigma)$ the set of 
all $P \in \P_K$ such that $F_Q = \sigma$ for some $Q \in \P_L$ 
lying above $P$. By Chebotarev's density theorem 
(see \cite[Theorem 7.30]{N}, \cite[Theorem 4.4.6]{HKalgnum}, or
\cite[Theorem 7.9.2]{HKclassfield}), the set
$\mathcal P(L/K,\sigma)$ has positive Dirichlet density. 
Therefore $\{L^{\langle \sigma \rangle} \mid \sigma \in G\}$ 
is a subset of the set of all decomposition fields, and hence
\[
K\subseteq K' \subseteq \bigcap_{\sigma \in G} L^{\langle \sigma \rangle} = K,
\]
which implies $K' = K$.

Now consider any $P\in\PK$ and some $Q\in\PL$ lying above $P$.
Let $Z$ be the decomposition field of $Q$ and $Q_Z = Q\cap \OZ$.
Then 
\[
\OKhash \subseteq K\cap \OZQZ = \OKP.
\]
Hence,
\[
\OK \subseteq \OKhash \subseteq \bigcap_{P\in\PK} \OKP = \OK,
\]
so that $\OKhash=\OK$.

If $L/K$ is an arbitrary finite extension, let $L^*/K$ be a finite 
Galois extension such that $L \subset L^*$. 

If $P \in \P_K$, \ $Q \in \P_L$ and $Q^* \in \P_{L^*}$ are such that 
$P \subset Q \subset Q^*$, let $C_Q$ be the $Q$-adic closure of 
$\OK$ in $L$ and $C_{Q^*}$ the $Q^*$-adic closure of $\OK$ in $L^*$.
Then $C_Q = C_{Q^*} \cap L$. Hence, as the intersection of all
$Q^*$-adic closures already equals $\OK$, all the more this holds 
for the intersection of all $Q$-adic closures.   
\end{proof}

\section{An application to integer-valued polynomials}

For an integral domain $D$ with quotient field $K$,
the ``ring of integer-valued polynomials over $D$'' consists of 
the polynomials with coefficients in $K$ that map elements of $D$
(when substituted for the variable) to elements of $D$:
\[
\Int(D) = \{f\in K[x] \mid f(D)\subseteq D \}.
\]

$P$-adic closure turns out to be useful for describing the image of 
an element of a number field under the ring of integer-valued 
polynomials of the ring of algebraic integers in a subfield. 
(A different description of the same image has been given by
McQuillan~\cite{McQuSplitPrimes93}).

Before showing this, we recall in gory detail, by request of the referee,
a standard argument for evaluating $v(f(c))$ when $c$ is an element of a 
discrete valuation ring $D$ and $f$ a product of monic linear factors in 
$D[x]$:

\begin{remark}\label{Legendre}
Consider Legendre's formula for the exponent of a prime $p$ in the
prime factorization of $n!$:
\[
\vp(1\cdot 2\cdot\ldots\cdot n) =
\sum_{k\ge 1}\biggl\lfloor\frac{n}{p^k}\biggr\rfloor.
\]

A priori, $\vp(1\cdot 2\cdot\ldots\cdot n) =\sum_{j=1}^n \vp(j)$,
but alternatively we can add up, for $k\ge 1$, the number of those
$j$ in $\{1,\ldots, n\}$ that are divisible by $p^k$, because in 
that way every $j$ with $\vp(j) = m$ is counted exactly $m$ times.

By the same token, for elements $c, a_1,\ldots, a_n$ in a discrete
valuation ring with maximal ideal $M$, valuation $v$, and valuation 
group $e\intz$, $e>0$,
\[
v(\prod_{j=1}^n (c-a_j)) = \sum_{j=1}^n v(c-a_j) =
e \sum_{k\ge 1}
\left|\{ 1\le j\le n \mid a_j\in c+M^k \} \right|.
\]
\end{remark}

\begin{proposition}\label{ImcinEQ}
Let $D$ be a Dedekind domain with finite residue fields,
and $K$ its quotient field. 

Let $c$ algebraic over $K$, $F=K[c]$ and 
$E$ the integral closure of $D$ in $F$. 

Then, for any maximal ideal $Q$ of $E$, 
the image of $c$ under $\Int(D)$ is contained in $E_Q$ 
if and only if $c$ is in the $Q$-adic closure of $D$ in $F$.
\end{proposition}

\begin{proof}
Let $P=Q\cap D$, and $PE= Q^e Q_1^{e_1}\ldots Q_{r-1}^{e_{r-1}}$
the prime factorization of $P$ in $E$.
Let $\vQ$ the valuation on $F$ associated to $Q$, normalized so 
that its value group is $\intz$, and, consequently, $e\intz$
the value group of its restriction to $K$, which we write as $\vP$.

Suppose $c$ is in the $Q$-adic closure of $D$ in $F$, and
$f\in\Int(D)$. Write $f=g/d$ with $g\in D[x]$ and $d\in D$.
Let $m=\vQ(d)$ and let $c'\in D$ such that $\vQ(c-c')\ge m$. 
Since $\vQ(g(c'))\ge m$ and $g(c')\equiv g(c)$ modulo $Q^m$,
we see $\vQ(g(c))\ge m$ and hence $f(c)\in E_Q$.

Conversely, suppose that $c$ is not in the $Q$-adic closure of $D$ in $F$.
Let $m\in \natn$ such that $(c + Q^{em})\cap D =\emptyset$. Let
$[D\colon P] =p$ and $a_1$, $\ldots,$ $a_{p^m}$ a complete system
of residues of $D$ modulo $P^m$.
Let $\beta = (1 - p^m)/(1-p)$, and $d\in K$ such that $\vP(d)=-e\beta$ 
and $v(d)\ge 0$ for all other essential valuations of $D$. 

Set $g(x)=\prod_{j=1}^{p^{m}}(x-a_j)$, and $f=dg$.
Then, by Remark~\ref{Legendre},
$\min_{r\in D}\vP(g(r))= e(1 + p +\ldots+ p^{(m-1)}) = e\beta$, the
minimum being attained by those $r$ with $\vP(r-a_j) = m$ for the 
unique $j$ such that $r\congr  a_j$ modulo $P^m$. 
Also, $f\in D_{P'}[x]$ for all maximal ideals $P'\ne P$ of $D$, whence
$f\in\Int(D)$. 

At the same time, $f(c)\notin E_Q$ since $\vQ(g(c))<e\beta$: 
To see this, calculate $\vQ(g(c))$ according to Remark~\ref{Legendre}:
If $c\notin E_Q$ then already $\vQ(g(c))<0$. So, assume $c\in E_Q$. then
$\vQ(g(c)) = \sum_{k\ge 1} n_k(c)$, where
\[
n_k(c) = \left|\{ 1\le j\le p^m \mid \vQ(c-a_j)\ge k \} \right|
= \left|\{ 1\le j\le p^m \mid a_j\in c+Q^k \} \right|.
\]
 Now, the intersection of a residue 
class of $Q^k$ in $E$ with $D$ is either empty or a residue class of 
$P^{\left\lceil\frac{k}{e}\right\rceil}$ in $D$, and $(c+Q^k)\cap D$ is 
empty for all $k\ge em$ by assumption. 

For $1\le s\le m$, each residue class of $P^s$ is represented among 
the $a_j$ exactly $p^{m-s}$ times, so that 
\[
\vQ(g(c)) = \sum_{k= 1}^{em-1} n_k(c) =
\sum_{s = 1}^{m-1} \sum_{r= 1}^{e} n_{(s-1)e+r}(c) + 
\sum_{r= 1}^{e-1} n_{(m-1)e+r}(c)
\]
implies
\[
\vQ(g(c)) \le
\sum_{s = 1}^{m-1}\sum_{r= 1}^{e} p^{m-s} + \sum_{r= 1}^{e-1} 1 =
\sum_{s = 1}^{m-1}e p^{m-s} + e - 1 = e\beta - 1
\]
\qedhere
\end{proof}

\begin{corollary}\label{imageintersection}
Let $D=\OK$ be the ring of integers in a number field $K$,
$c$ algebraic over $K$, $F=K[c]$, and
$E=\OF$ the integral closure of $D$ in $F$.

Then the image of $c$ under $\Int(D)$ is 
\[
\Int(D)[c] = \bigcap_{Q\in \Pc} E_Q,
\]
where $\Pc$ is the set of those $Q\in\Spec(E)$ such that $c$ is 
in the $Q$-adic closure of $D$.
\end{corollary}

\begin{proof}
Clearly, $D[c]\subseteq \Int(D)[c]\subseteq K[c]=F$.
Since $\Int(D)$ is Pr\"ufer, so is $\Int(D)[c]$, as a homomorphic
image of a Pr\"ufer domain. In particular, $\Int(D)[c]$ is
integrally closed in its quotient field $F$, and, therefore,
contains $E$. As an overring of the Dedekind ring $E$, $\Int(D)[c]$ 
is necessarily an intersection of localizations of $E$ at maximal
ideals, and, therefore, equal to the intersection of all $E_Q$
containing it. Proposition~\ref{ImcinEQ} tells us which ones these are.
\end{proof}

\begin{corollary}\label{imageoverring}
Let $D=\OK$ be the ring of integers in a number field $K$,
$F$ a finite extension of $K$, and $E=\OF$ the integral
closure of $D$ in $F$. Let $c\in F$.

Then $\Int(D)[c]$, the image of $c$ under $\Int(D)$, is
\[
\Int(D)[c]=
\begin{cases}
D&{\text{if}}\quad c\in D \\
{\text{a strict overring of $D$}}&{\text{if}}\quad c\in K\setminus D \\
{\text{a strict overring of $E$}}&{\text{if}}\quad c\in F\setminus K \\
\end{cases}
\]
\end{corollary}

\begin{proof}
The first two cases are trivial; the third follows from 
Theorem~\ref{intersectionOK} and Corollary~\ref{imageintersection}.
\end{proof}

The relevance of the above corollary is this: it allows us, by applying 
a result of Vaserstein \cite{VsteinDedekEN}, to conclude that the 
residue class rings of $\IntOK$ modulo non-zero prime ideals are
$\GEt$, for those rings of integers $\OK$ that are
themselves $\GEt$. (We are interested in showing this
as residue class rings that are $\GEt$ can be used as
a step towards determining the stable rank of $\IntOK$.)

Generalized Euclidean rings were introduced by Cohn in a seminal paper
\cite{Cohn66GL2} in 1966. A commutative ring $R$ is $\GEt$ if any
unimodular pair $(a,b)\in R^2$ (that is, any pair such that $aR + bR = R$) 
can be transformed to $(1,0)$ by a series of elementary transformations, 
where an elementary transformation consists of replacing $(a,b)$ by
$(a,b+ra)$ or by $(a+rb,b)$ for some $r\in R$. Likewise, $R$ is 
called $\GEn$ if every unimodular $n$-tuple can be transformed to
$(1,0,\ldots,0)$ by a series of elementary transformations (consisting
of adding a scalar multiple of one entry to a different entry), and $R$
is called generalized Euclidean if it is $\GEn$ for all $n>0$.
By the Euclidean algorithm, Euclidean rings are generalized Euclidean.

Vaserstein showed for the ring of integers $\OK$ in a number field
$K$ that firstly, every strict overring of $\OK$ is $\GEt$, and, 
secondly, when $K$ is not imaginary quadratic then $\OK$ itself is 
$\GEt$.

Cohn~\cite{Cohn66GL2} had already shown that in an imaginary quadratic 
number field $K$, the ring of integers $\OK$ is not $\GEt$ unless it is
actually Euclidean. Since the Euclidean imaginary quadratic $\OK$ are 
known, we have in the results of Cohn~\cite{Cohn66GL2} and 
Vaserstein~\cite{VsteinDedekEN} a complete classification of those
rings of integers in number fields (and their overrings) that are 
$\GEt$.

\begin{corollary}
Let $K$ be one of those number fields for which $\OK$ is $\GEt$. Then 
$\IntOK/P$ is $\GEt$ for every non-zero prime ideal $P$ of $\IntOK$.

In particular, $\IntOK/P$ is $\GEt$ for every non-zero prime ideal $P$ 
of $\IntOK$ whenever $\OK$ is Euclidean or $K$ not imaginary quadratic.
\end{corollary}

\begin{proof}
If $P$ is maximal then $\IntOK/P$ is a field, and hence Euclidean.
Any non-zero non-maximal prime ideal of $\IntOK$ is of the form
$\IntOK \cap f(x)K[x]$ for an irreducible polynomial $f\in K[x]$.
Let $c$ be a root of $f$ in the splitting field of $f$ over $K$.
Clearly then $\IntOK/P$ is isomorphic to the image of $c$ under $\IntOK$,
which is $\GEt$ by Corollary~\ref{imageoverring} and the hypothesis.
\end{proof}

{\textbf{Acknowledgment:} }
The surviving author wishes to thank the referee for graciously 
correcting some inaccuracies due to Prof.~Halter-Koch's sudden 
death in the drafting stage of the paper, before proofreading.

\frenchspacing
% bibliography
\bibliographystyle{siamese}
\bibliography{ivimagebib}
\nocite{*}

\end{document}